\documentclass{article}
\usepackage{graphicx}
\def\proc#1{\medbreak\noindent{\bf #1}\ \ignorespaces}
\def\ep{\noindent{\hfill$\diamondsuit$}}
\def\address#1{\\ #1}
\def\email#1{\footnote{email: #1}}
\def\recd#1{\date{#1}}
\newtheorem {theorem} {Theorem}
\newtheorem {prop}{Proposition}[section]
\newtheorem {lemma}[prop]{Lemma}
\newtheorem {corollary}[prop]{Corollary}
\newtheorem {definition}[prop]{Definition} 
\def\baselinestretch{1.2}
\def\baselinestretch{1.2}
\begin{document}

\title{ Elliptic Islands on Strictly Convex Billiards}

\author{M\'ario Jorge DIAS CARNEIRO, Sylvie OLIFFSON KAMPHORST \\
and S\^onia PINTO DE CARVALHO
\address{Departamento de Matem\'atica, ICEx, UFMG \break
 30.123--970, Belo Horizonte, Brasil.}
\email{carneiro@mat.ufmg.br, syok@mat.ufmg.br, sonia@mat.ufmg.br}
}

\recd{January 2002}

\maketitle

\begin{abstract}
This paper addresses the question of genericity of existence of elliptic
islands for the billiard map associated to strictly convex closed curves. 
More precisely, we study 2-periodic orbits of billiards associated to $C^5$
closed and strictly convex curves and show that the existence of elliptic
islands is a dense property on the subset of those billiards having an
elliptic 2-periodic point.

Our main tools are normal perturbations, the Birkhoff Normal Form for elliptic
fixed points and Moser's Twist Theorem.
\end{abstract}

\section{Introduction}

Let $\alpha$ be a plane, closed, regular,
strictly convex, oriented counterclockwise.
Any curve with these properties may be parametrized by
$\varphi\in [0,2\pi]$, the angle between the tangent vector and a given
direction
The billiard problem on $\alpha$ consists in the free motion of a point
particle in the plane region enclosed by $\alpha$, being reflected
elastically at the impacts with the boundary.
Since the particle moves with constant velocity inside
the region, the motion is completely determined  by the point of  reflection at
$\alpha$ and the
direction of movement immediately after each reflection. 
So, 
the angle $\varphi$, which locates the point of  reflection, and the
angle $\theta$ between the direction of motion and the normal to the boundary
at the reflection point, may be used to describe the system.

The billiard model defines a map $T$ from the annulus 
${\cal A}=[0,2\pi)\times(-\pi/2,\pi/2)$ into itself. 
If $\alpha$ is $C^k$, $T$ is a $C^{k-1}$-diffeomorphism \cite{kn:str} preserving
the
measure  $d\mu = R(\varphi) \cos\theta \, d\theta d\varphi$ \cite{kn:bir}, 
where
$R(\varphi)$ is the radius of curvature of $\alpha$ at $\varphi$. 
$({\cal A}, \mu, T)$ defines a two dimensional discrete dynamical system,
whose orbits are  given by  $\{ (\varphi_n,
\theta _n)= T^n(\varphi_0,\theta_0), n \in Z  \}. $

A point $(\varphi_0,\theta_0)\in {\cal A}$ is $n$-periodic if $n$ is the
smallest positive integer such that $T^n(\varphi_0,\theta_0)=
(\varphi_0,\theta_0)$. It is called elliptic if the eigenvalues of
$DT^n_{(\varphi_0,\theta_0)}$ are complex, with non vanishing real part. An
elliptic island is a $T^n$-invariant subset, homeomorphic to a disk,
surrounding a $n$-periodic elliptic point.

We study 2-periodic orbits of billiards associated to $C^5$
closed and strictly convex curves and show that the existence of elliptic
islands is a dense property on the subset of those billiards having an
elliptic 2-periodic point.

It is well known that for strictly convex curves, at least $C^6$, the
billiard map has invariant curves near the boundary of the phase space
\cite{kn:dou}, and then it is not ergodic. Here we study the dynamical
properties of the region of the phase space between these curves, when they exist,
or the whole phase space, otherwise.
The existence of islands assures the non-ergodicity of this region.

The layout of the paper is as follows: 
we will deal with plane, closed, regular,
strictly convex, oriented counterclockwise curves.
In sections 2 and 3 we analyze the geometrical
and dynamical properties of the 2-periodic orbits. In section 4, we define
normal perturbations and prove that the existence of elliptic or hyperbolic
2-periodic orbits is an open property. Sections  5 and 6 deal with the resonances
and the existence of elliptic islands. In sections 7 we present a more topological
description of the results obtained in the previous sections. We conclude this
work with some comments and open questions.
Our main tools are normal perturbations, the Birkhoff Normal Form for elliptic
fixed points and Moser's Twist Theorem.

\section{The Diameters of a Convex Curve}

Let a curve $\alpha$ be parametrized by
$\varphi\in [0,2\pi]$, the angle between the tangent vector and a given
direction, hereafter denoted by $\bf x$.
Let
$\tau(\varphi)$ and $\eta(\varphi)$ be the unitary tangent and
(inward) normal vectors. Then $\alpha'(\varphi) = R(\varphi) \tau (\varphi)$. 
Clearly,
$\tau(\varphi+\pi)=-\tau(\varphi)$ and 
$\eta(\varphi+\pi)=-\eta(\varphi)$.

It is known that any closed convex curve has at least a width and a diameter,
characterized by the global maximum and minimum of the distance between points
with parallel tangents (or with double normal) (see, for instance,
\cite{kn:kat}). We broaden this concept defining:
\begin{definition}  
$\varphi_0$ is a diameter of
$\alpha$ if $\alpha(\varphi_0)-\alpha(\varphi_0+\pi)$ is parallel to
$\eta(\varphi_0)$.  
\end{definition}

\begin{prop}\label{prop:diam}
Let $l(\varphi)=|\alpha(\varphi)-\alpha(\varphi+\pi)|$. Then $\varphi_0$ is a
diameter of $\alpha$ if and only if
$\frac{dl}{d\varphi}(\varphi_0)=0$
\end{prop}

\proc{Proof:}
\begin{eqnarray*}
2l \frac{dl}{d\varphi} &=&
2  \, \langle \alpha'(\varphi)-\alpha'(\varphi+\pi), 
\alpha(\varphi)-\alpha(\varphi+\pi)
\rangle = \\
&=& 2 \, \langle R(\varphi)\tau(\varphi)-R(\varphi+\pi)\tau(\varphi+\pi),
\alpha(\varphi)-\alpha(\varphi+\pi) \rangle = \\
&=& 2 \, \langle (R(\varphi)+R(\varphi+\pi))\tau(\varphi),\alpha(\varphi)-
\alpha(\varphi+\pi) \rangle
\end{eqnarray*}
So $$\frac{dl}{d\varphi}(\varphi_0)=0\Longleftrightarrow 
\langle \alpha(\varphi_0)-\alpha(\varphi_0+\pi),\tau(\varphi_0) \rangle = 0$$ 
\ep\medbreak

Let $L_0=l(\varphi_0)$, $R_0=R(\varphi_0)$ and
$R_\pi=R(\varphi_0+\pi)$.
\begin{prop}
Let $\varphi_0$ be a diameter.
If $L_0-(R_0+R_\pi)\neq 0$ then $\varphi_0$ is an isolated diameter.
\end{prop}

\proc{Proof:} 
\begin{eqnarray*}
&& \frac{d}{d\varphi}\left(2l\frac{dl}{d\varphi}\right)
 =  2 \left( \frac{dl}{d\varphi} \right)^2 + 2l \frac{d^2l}{d\varphi^2} = \\
&&\   =  2 \, ( \langle (R'(\varphi)+R'(\varphi+\pi))\tau(\varphi)+
(R(\varphi)+R(\varphi+\pi))
\eta(\varphi),\alpha(\varphi)-\alpha(\varphi+\pi) \rangle +\\
&& \ \ \ \ + \langle ( R(\varphi)+R(\varphi+\pi))\tau(\varphi),(R(\varphi)+
R(\varphi+\pi))\tau(\varphi) \rangle )
\end{eqnarray*}
At the diameter $\varphi_0$:
$$2l(\varphi_0)\frac{d^2l}{d\varphi^2}(\varphi_0)=2(-L_0(R_0+R_\pi)+(R_0+R_\pi)^2)=
-2(R_0+R_\pi)(L_0-(R_0+R_\pi)).$$
and we have that if $L_0-(R_0+R_\pi)\neq 0$ then $\varphi_0$ is an isolated
singular point of $l$.
\ep\medbreak

\section{Dynamical Characterization of the Diameters of a Convex Curve}

Let $\alpha$ be a $C^k$ curve, $k\geq 2$,
and $T$ its associated billiard map.

If $\varphi_0$ is a diameter for $\alpha$ then $(\varphi_0,0)$ is a 2-periodic
orbit for $T$, and, as we have pointed out, any billiard on a strictly convex
curve has at least two of those orbits. Its stability is given by
the eigenvalues of: 
$$
DT^2_{(\varphi_0,0)}=\frac{1}{R_0R_\pi}\left( 
\begin{array}{cc} 
L_0-R_\pi&L_0\\ 
L_0-(R_0+R_\pi)& L_0-R_0
\end{array}\right)
\left(  \begin{array}{cc} 
L_0-R_0&L_0\\ 
L_0-(R_0+R_\pi)&L_0-R_\pi
\end{array}\right) 
$$
it follows that:
\begin{enumerate}
\item $(\varphi_0,0)$ is hyperbolic if
\begin{itemize}
\item $L_0-(R_0+R_\pi)>0$ or
\item $L_0-(R_0+R_\pi)<0$ and $(L_0-R_0)(L_0-R_\pi)<0$.
\end{itemize}
\item $(\varphi_0,0)$ is elliptic if $L_0-(R_0+R_\pi)<0$ 
and $(L_0-R_0)(L_0-R_\pi)>0$.
\item $(\varphi_0,0)$ is parabolic if
\begin{itemize}
\item $L_0-(R_0+R_\pi)=0$ or
\item $L_0-(R_0+R_\pi)<0$ and $(L_0-R_0)(L_0-R_\pi)=0$.
\end{itemize}
\end{enumerate}
We will say that a diameter is hyperbolic, elliptic or parabolic if the
associated 2-periodic orbit is, respectively, hyperbolic, elliptic or
parabolic. 

A nice consequence \cite{kn:koz} of this characterization is that, for $k\geq
3$, the largest diameter of a strictly convex curve, if isolated, is
hyperbolic, due to the fact that
$
\displaystyle\frac{d^2l}{d\varphi^2}(\varphi_0)
=-\frac{(R_0+R_\pi)(L_0-(R_0+R_\pi))}{L_0} < 0$.

The smallest one, however, can be elliptic or hyperbolic. If $R_0=R_\pi$ (which
happens, for instance, for the ellipse), it is elliptic.
Nevertheless, as we show in the last section, there are strictly convex curves
such that all the diameters are hyperbolic.

When $\varphi_0$ is elliptic, the eigenvalues of $DT^2_{(\varphi_0,0)}$ are given
by $\displaystyle e^{\pm i\gamma}=\frac{\mbox{tr}(DT^2)}{2}\pm
i\sqrt{1-(\frac{\mbox{tr}(DT^2)}{2})^2}$ and so
\begin{equation}\label{eq:autov}
\cos\gamma=2\frac{(L_0-R_0)(L_0-R_\pi)}{R_0R_\pi}-1 \ .
\end{equation}

\section {Normal Perturbations of a Convex Curve}

Let $\alpha$ 
$C^k$, with $k\geq 4$, when parametrized by
$\varphi\in [0,2\pi]$, the angle between the tangent vector and $\bf x$.
Using this direction and its perpendicular as a reference frame, 
the unitary tangent and (inward) normal vectors are written, respectively, as 
$\tau(\varphi)=(\cos\varphi,\sin\varphi)$ and
$\eta(\varphi)=(-\sin\varphi,\cos\varphi)$. So, the normal bundle
$(\alpha(\varphi),\eta(\varphi))$ is also $C^k$ and $\eta(\varphi)$ is
$C^\infty$. 

A normal perturbation of $\alpha$ is a curve $\beta$ given by 
$\beta(\varphi)=\alpha(\varphi)+\lambda(\varphi)\eta(\varphi)$, where 
$\lambda$ is a $C^k$, $2\pi$-periodic function. Observe that in general $\varphi$
does not represent the angle between the tangent vector of $\beta$ and the
fixed direction $\bf x$.

We will use the norm  
$ \displaystyle |\lambda|_2 = \max_{\varphi\in
[0,2\pi]}\{|\lambda(\varphi)|,|\lambda'(\varphi)|,|\lambda''(\varphi)|\}$. 

\begin{prop}
If  $|\lambda|_2$ is small enough then $\beta$
is  $C^k$, closed, regular and strictly convex.
\end{prop}

\proc{Proof:} : $\beta$ is closed because $\lambda$ is $2\pi$-periodic.
$\alpha$ is $C^k$, $\lambda$ is $C^k$ and $\eta$ is $C^\infty$. So, $\beta$ is
$C^k$. $\beta '=\alpha '+\lambda '\eta-\lambda\tau$. So, if $\alpha$ is
regular, for $\epsilon$ small enough, $\beta$ is regular.
$\beta ''\times\beta '=(\alpha ''\times\alpha ')+\lambda '(\alpha
''\times\eta)-\lambda(\alpha ''\times\tau)+(\lambda
''-\lambda)(\eta\times\alpha ')+(\lambda(\lambda ''-\lambda)+2\lambda
')(\tau\times\eta)$. So, if $\alpha$ is strictly convex, for
$\epsilon$ small enough, $\beta$ is strictly convex. 
\ep\medbreak

\begin{figure}
\centering
\includegraphics[width=7.5cm,height=4cm]{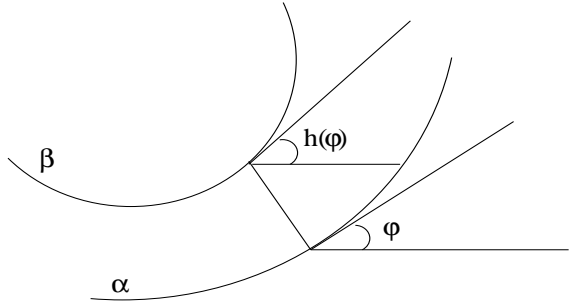}
\caption{Definition of $h(\varphi)$}
\label{fig:h}
\end{figure}

Let us call $\vec{t}(\varphi)$ and $\vec{n}(\varphi)$ the unitary tangent and
(inward) normal vectors and $R_{\alpha}(\varphi)$ and $R_{\beta}(\varphi)$ the
radius of curvature of $\alpha$ and $\beta$ at $\varphi$, respectively. Then:
\begin{eqnarray}
\vec{t} &=& \frac{(R_{\alpha}-\lambda)\tau+\lambda '\eta}
{((R_{\alpha}-\lambda)^2+(\lambda ')^2)^{1/2}} \nonumber\\
\vec{n} &=& \frac{-\lambda '\tau+(R_{\alpha}-\lambda)\eta}
{((R_{\alpha}-\lambda)^2+(\lambda ')^2)^{1/2}} \nonumber\\
\label{eq:rbeta}
R_{\beta} &=& \frac{((R_{\alpha}-\lambda)^2+(\lambda ')^2)^{3/2}}
{(R_{\alpha}-\lambda)(R_{\alpha}-\lambda+\lambda '')-\lambda '
(R'_{\alpha}-2\lambda ')}
\end{eqnarray}
Let $h(\varphi)$ be the angle between $\vec{t}(\varphi)$ and
the given direction that defines $\varphi$, as in fig.~\ref{fig:h}.
It is easy to see that 
$$h(\varphi)=\varphi+\arctan\frac{\lambda '}{R_{\alpha}-\lambda}.$$

According to proposition \ref{prop:diam}, the diameters of $\beta$
are such that the derivative of
$$l_{\beta}(\varphi)=|\beta(\varphi)-\beta(h^{-1}(h(\varphi)+\pi))|$$
vanishes.

\begin{lemma} \label{lema:difeo}
If $|\lambda|_2$ is sufficiently small
then $h$ is a $C^2$-diffeomorphism.  
\end{lemma}

\proc{Proof:} : As $\alpha$ is closed, $0<c_1\leq R_{\alpha}$. Then, for $|\lambda|_2$
sufficiently small, $R_{\alpha}-\lambda>0$ and $h$ is well defined. 
$\alpha$ and $\beta$ being at least $C^4$ then $\lambda$ is at least $C^4$,
and so $h$ is at least $C^3$.

Let $s$ be the arclength parameter for $\beta$. Then
$$\frac{dh}{d\varphi}=
\frac{dh}{ds}\frac{ds}{d\varphi}=k_{\beta}|\frac{d\beta}{d\varphi}|\geq
c_2>0.$$ Then $h$ is $C^1$ and invertible, and $h^{-1}$ has continuous
derivative.
$$\frac{d^2h}{d\varphi^2}=\frac{dk_{\beta}}{ds}|\frac{d\beta}{d\varphi}|^2
+k_{\beta}\frac{<\frac{d^2\beta}{d\varphi^2},\frac{d\beta}{d\varphi}>}
{|\frac{d\beta}{d\varphi}|}$$
so $h$ is $C^2$.
Taking $y=h(\varphi)$ we have
$$\frac{d^2h^{-1}}{dy^2}(y)=-\frac{\frac{d^2h}{d\varphi^2}(\varphi)}
{(\frac{dh}{d\varphi}(\varphi))^3}$$
which shows that $h$ is a $C^2$-diffeomorphism.
\ep\medbreak

\begin{lemma} \label{lem:h-1} 
With the same hypothesis of lemma
\ref{lema:difeo},   $f(\varphi)=h^{-1}(h(\varphi)+\pi)$ and
$I_{\pi}(\varphi)=\varphi+\pi$ are $C^2$ close.
\end{lemma}

\proc{Proof:} : Let $|\lambda|_2 < \epsilon$ and
$\frac{1}{c_1}<|R_{\alpha}-\lambda|<c_2$,  then
$|h(\varphi)-\varphi|=|\arctan\frac{\lambda '}{R_{\alpha}-\lambda}|\leq
|\frac{\lambda '}{R_{\alpha}-\lambda}|<c_1\epsilon$

Let  $c_3=\hbox{max}|\frac{dh^{-1}}{dy}|$.
$$ |f(\varphi)-I_{\pi}(\varphi)|=|h^{-1}(h(\varphi)+\pi)-h^{-1}(h(\varphi+\pi))|\leq
c_3|h(\varphi)+\pi-h(\varphi+\pi)|=$$
$$=c_3|h(\varphi)-\varphi+\varphi+\pi-h(\varphi+\pi)|
\leq c_3(|h(\varphi)-\varphi|+|\varphi+\pi-h(\varphi+\pi)|)<2c_1c_3\epsilon
$$
Let
$0<\frac{1}{c_4}<|R'_{\alpha}-\lambda '|<c_5$ and 
$\Delta(\varphi)=\arctan\frac{\lambda '}{R_{\alpha}-\lambda}$.
Now, $\frac{dh}{d\varphi}=1+\frac{d\Delta}{d\varphi}$ and
\begin{equation}\label{eq:del}
|\frac{d\Delta}{d\varphi}|\leq|\frac{d}{d\varphi}\frac{\lambda
'}{R_{\alpha}-\lambda}|=|\frac{\lambda ''(R_{\alpha}-\lambda)-\lambda
'(R'_{\alpha}-\lambda ')}{(R_{\alpha}-\lambda)^2}|< c_1^2c_2c_5\epsilon
\end{equation}
Taking $y=h(\varphi)$, 
\begin{eqnarray*}
\left| \frac{df}{d\varphi}-1 \right|&=&
\left| \frac{dh^{-1}}{dy}(h(\varphi)+\pi)\frac{dh}{d\varphi}(\varphi)-1 \right|
=\left| \frac{1+\frac{d\Delta}{d\varphi}(\varphi)}
{1+\frac{d\Delta}{d\varphi}(h^{-1}(h(\varphi)+\pi))}-1 \right|=\\
&=&\left| \frac{1+\frac{d\Delta}{d\varphi}(\varphi)
+\frac{d\Delta}{d\varphi}(h^{-1}(h(\varphi)+\pi))
-\frac{d\Delta}{d\varphi}(h^{-1}(h(\varphi)+\pi))}
{1+\frac{d\Delta}{d\varphi}(h^{-1}(h(\varphi)+\pi))}-1 \right| = \\
&=& \left| \frac{\frac{d\Delta}{d\varphi}(\varphi)
-\frac{d\Delta}{d\varphi}(h^{-1}(h(\varphi)+\pi))}
{1+\frac{d\Delta}{d\varphi}(h^{-1}(h(\varphi)+\pi))} \right| \leq\\
&\leq& \left| \frac{dh^{-1}}{d\varphi}(h(\varphi)+\pi)) \right|
\, \left| \frac{d\Delta}{d\varphi}(\varphi) \right| +
\left| \frac{d\Delta}{d\varphi}(h^{-1}(h(\varphi)+\pi)) \right|
<2c_3c_1^2c_2c_5\epsilon
\end{eqnarray*}
Finally,
\begin{eqnarray*}
|\frac{d^2f}{d\varphi^2}(\varphi)-0|=|\frac{d}{d\varphi}
(\frac{dh^{-1}}{dy}(h(\varphi)+\pi)\frac{dh}{d\varphi}(\varphi))|\leq
\\ |\frac{d^2h^{-1}}{dy^2}(h(\varphi)+\pi)||\frac{dh}{d\varphi}|^2+
|\frac{dh^{-1}}{dy}(h(\varphi)+\pi)||\frac{d^2h}{d\varphi^2}(\varphi))|
\end{eqnarray*}
$|\frac{dh}{d\varphi}|$ and $|\frac{dh^{-1}}{dy}(h(\varphi)+\pi)|$ are bounded
from above and bellow because $h$ is a $C^2$ diffeomorphism on $[0,2\pi]$.
\begin{eqnarray}
&|\frac{d^2h}{d\varphi^2}(\varphi))|=|\frac{d^2\Delta}{d\varphi^2}(\varphi))|\leq
\nonumber\\
&|\frac{\lambda '''(R_{\alpha}-\lambda)-\lambda '(R''_{\alpha}-\lambda '')}
{(R_{\alpha}-\lambda)^2+(\lambda ')^2}|
+2|\frac{(\lambda ''(R_{\alpha}-\lambda)-\lambda
'(R'_{\alpha}-\lambda '))(R'_{\alpha}-\lambda ')}
{(R_{\alpha}-\lambda)((R_{\alpha}-\lambda)^2+(\lambda ')^2)}|
+2|\frac{\lambda
'}{R_{\alpha}-\lambda}||\frac{d\Delta}{d\varphi}(\varphi)|\nonumber
\end{eqnarray} 
For $\epsilon$ small enough, $0<\frac{1}{c_1}<|R_{\alpha}-\lambda|<c_2$, 
 $|R'_{\alpha}-\lambda '|<c_5$, $|R''_{\alpha}-\lambda ''|<c_6$ and 
$(R_{\alpha}-\lambda)^2+(\lambda ')^2=|\beta '(\varphi)|^2\geq c_7>0$, since
$\beta$ is regular. Together with equation \ref{eq:del}, we have that
$|\frac{d^2h}{d\varphi^2}(\varphi))|<K\epsilon$, for a constant $K$.
\begin{eqnarray}
|\frac{d^2h^{-1}}{dy^2}(h(\varphi)+\pi)|=
|\frac{\frac{d^2h}{d\varphi^2}(h^{-1}(h(\varphi)+\pi))}
{(\frac{dh}{d\varphi}(h^{-1}(h(\varphi)+\pi)))^3}|<K\epsilon
\end{eqnarray}
for a constant $K$.
So, $|\frac{d^2f}{d\varphi^2}(\varphi)|<C\epsilon$.
\ep\medbreak

\begin{lemma}\label{lem:lalfa}
 Let
$l_{\alpha}(\varphi)=|\alpha(\varphi)-\alpha(\varphi+\pi)|$ and 
$l_{\beta}(\varphi)=|\beta(\varphi)-\beta(h^{-1}(h(\varphi)+\pi))|$.
If $|\lambda|_2$ is small enough
then $l_{\alpha}$ and $l_{\beta}$ are $C^2$ close.
\end{lemma}

\proc{Proof:} : If $|\lambda|_2$ is sufficiently small then $\alpha$ and $\beta$ are
$C^2$ close. So $\alpha\circ I_{\pi}$ and $\beta\circ f$ are $C^2$ close and
the result follows.
\ep\medbreak

\begin{prop} \label{prop:iso}
Suppose that $\alpha$ has an isolated diameter at $\varphi_0$, with
$\frac{dl_{\alpha}}{d\varphi}(\varphi_0)=0$ and
$\frac{d^2l_{\alpha}}{d\varphi^2}(\varphi_0)>0$ (or $<0$). Let $\delta$ be such
that $\frac{d^2l_{\alpha}}{d\varphi^2}(\varphi)>0$ (or $<0$) for $\varphi\in
[\varphi_0-\delta,\varphi_0+\delta]$. Then there exists $\epsilon_0$ such that
if $|\lambda|_2<\epsilon_0$ then any normal perturbation 
$\beta(\varphi)=\alpha(\varphi)+\lambda(\varphi)\eta(\varphi)$ has an
isolated diameter $\varphi_1 \in (\varphi_0-\delta,\varphi_0+\delta)$, with 
$\frac{d^2l_{\beta}}{d\varphi^2}(\varphi_1)>0$ (or $<0$).
Furthermore, $\epsilon_0\to 0$ if $\delta\to 0$. 
\end{prop}

\proc{Proof:} : The proof follows immediately from the results above and from the
transversality of the zero of $\frac{dl_{\alpha}}{d\varphi}$. 
\ep\medbreak

Observe that rotating the axis, if necessary, we can always suppose that the
diameter of $\alpha$ is at $\varphi=0$.

Let, as above, $L_0=l_{\alpha}(0)$, $R_0=R_{\alpha}(0)$,
and $R_{\pi}=R_{\alpha}(\pi)$. An obvious consequence is:
\begin{corollary}\label{coro:hip}
If $0$ is a hyperbolic diameter of $\alpha$, with $L_0-(R_0+R_\pi)>0$
(meaning that $\frac{d^2l_{\alpha}}{d\varphi^2}(0)<0$) then every normal
perturbation $\beta$ with $|\lambda|_2<\epsilon_0$ has a hyperbolic diameter
close to $0$.
\end{corollary}

\begin{prop}\label{prop:aberto}
Let $\epsilon_0$ and $\delta$ be as in Proposition \ref{prop:iso} and suppose
$0$ is an elliptic (resp. hyperbolic) diameter for $\alpha$, i.e.
$L_0-(R_0+R_\pi)<0$ and $(L_0-R_0)(L_0-R_\pi)>0$ (resp.
$(L_0-R_0)(L_0-R_\pi)<0$).  There exists $\epsilon _1\leq \epsilon_0$ such
that if $|\lambda|_2<\epsilon_1$ then any normal perturbation 
$\beta(\varphi)=\alpha(\varphi)+\lambda(\varphi)\eta(\varphi)$ has an elliptic
(resp. hyperbolic) diameter at $\varphi_1 \in (\varphi_0-\delta,\varphi_0+\delta)$.  
\end{prop} 

\proc{Proof:}
If
$|\lambda|_2<\epsilon_0$ then, by Proposition \ref{prop:iso}, $\beta$ has an
isolated diameter at $\varphi_1\in (\varphi_0-\delta,\varphi_0+\delta)$ with
$l_\beta(\varphi_1)-(R_\beta(\varphi_1)+R_\beta(h^{-1}(h(\varphi_1)+\pi))<0$.
Since $\alpha$ and $\beta$ are $C^2$-close, $R_\alpha$ and $R_\beta$ are
$C^0$-close. This, together with lemma \ref{lem:h-1} and
$(L_0-R_0)(L_0-R_\pi)>0$, gives that there exists $\epsilon _1\leq \epsilon_0$
such that if $|\lambda|_2<\epsilon_1$ then
$(l_\beta(\varphi_1)-R_\beta(\varphi_1))
(l_\beta(\varphi_1)-R_\beta(h^{-1}(h(\varphi_1)+\pi)))>0$.
\ep\medbreak

The hyperbolic case $(L_0-R_0)(L_0-R_\pi)<0$ is analogous.
\proc{Remark:} Propositions \ref{prop:iso} and \ref{prop:aberto}
and corollary \ref{coro:hip} only express, in the context of billiard
maps, the well known fact of openness of existence of elliptic or
hyperbolic fixed points. 

\section{Resonances}

We will assume that the curve $\alpha$ has an elliptic diameter with 
eigenvalues $e^{\pm i\gamma}$. Without any loss of generality we can assume
that this diameter is located at $\varphi = 0$.

\begin{definition}\label{def:res}
An elliptic diameter is resonant if its eigenvalues satisfy
$(e^{\pm i\gamma})^j=1$, for $j=2,3,4$, i.e., if
$\gamma=\frac{1}{2}\pi$, $\frac{2}{3}\pi$, $\frac{4}{3}\pi$, 
$\frac{3}{2}\pi$.
Otherwise, it is non-resonant.
\end{definition}

\begin{prop}\label{prop:reso}
Given a $C^k$ strictly convex curve $\alpha$, $k\geq 4$, with a resonant
elliptic diameter $\varphi=0$,
there is a $C^k$ strictly convex curve $\beta$, $C^2$-close
to $\alpha$, such that $\beta$ has a non-resonant elliptic diameter. 
\end{prop}

\proc{Proof:}
For $\epsilon_1$ as in Proposition \ref{prop:aberto} consider
the normal perturbations
$\beta(\varphi)=\alpha(\varphi)+\lambda(\varphi)\eta(\varphi)$, with
$|\lambda|_2<\epsilon_1$ and $\lambda '(0)=\lambda '(\pi)=0$, which implies
that $\varphi=0$ is also an elliptic diameter for $\beta$. 

Let $\lambda(0)=\lambda_0$,
$\lambda(\pi)=\lambda_{\pi}$, $\lambda''(0)=\lambda''_0$, $\lambda
''(\pi)=\lambda ''_{\pi}$ . We have that
$l_0=l_{\beta}(0)=|\beta(0)-\beta(\pi)|=L_0-\lambda_0-\lambda_{\pi}$ and, by
equation (\ref {eq:rbeta}),
 $$r_0=R_{\beta}(0)=\frac{(R_0-\lambda_0)^2}
{R_0-\lambda_0+\lambda ''_0}= R_0-\lambda_0- \lambda ''_0\frac{R_0-\lambda_0}
{R_0-\lambda_0+\lambda ''_0}$$
and
$$ r_{\pi}=R_{\beta}(\pi)=\frac{(R_{\pi}-\lambda_{\pi})^2}
{R_{\pi}-\lambda_{\pi}+\lambda ''_{\pi}}=R_\pi-\lambda_\pi- \lambda
''_\pi\frac{R_\pi-\lambda_\pi} {R_\pi-\lambda_\pi+\lambda ''_\pi}
$$
So, for $|\lambda|_2$ small enough, these normal perturbations
have an elliptic diameter with 
eigenvalues $e^{\pm i\gamma(\beta)}$ satisfying, by equation (\ref{eq:autov})
\begin{eqnarray*}
\cos\gamma(\beta)&=&2\,\frac{(l_0-r_0)(l_0-r_\pi)}{r_0r_\pi}-1=\\
&=&
2\,\frac{(L_0-R_0-\lambda_{\pi}+\lambda_0''\Delta_0)(L_0-R_\pi-\lambda_0+\lambda
''_\pi\Delta_\pi)}{(R_0-\lambda_0- \lambda ''_0\Delta_0)(R_\pi-\lambda_\pi-
\lambda ''_\pi\Delta_\pi)}-1=\\
&=&F(\lambda_0,\lambda_\pi,\lambda_0',\lambda_\pi')
\end{eqnarray*}
where $$\Delta_a=\frac{R_a-\lambda_a}
{R_a-\lambda_a+\lambda ''_a} \ .
$$
Clearly $F(0,0,0,0)=2\frac{(L_0-R_0)(L_0-R_\pi)}{R_0R_\pi}-1$
gives the eigenvalues of the elliptic diameter of $\alpha$.

Moreover
\begin{eqnarray*}
&&\frac{\partial F}{\partial\lambda_0}(0,0,0,0)=
2\frac{(L_0-R_0)(L_0-R_0-R_\pi)}{R_0^2R_\pi}\neq 0
\\
&&\frac{\partial F}{\partial\lambda_\pi}(0,0,0,0)=
2\frac{(L_0-R_\pi)(L_0-R_0-R_\pi)}{R_\pi^2R_0}\neq 0
\\
&&\frac{\partial F}{\partial\lambda_0''}(0,0,0,0)=
2\frac{(L_0-R_\pi)L_0}{R_0^2R_\pi}\neq 0
\\ 
&&\frac{\partial F}{\partial\lambda_\pi''}(0,0,0,0)=
2\frac{(L_0-R_0)L_0}{R_\pi^2R_0}\neq 0
\end{eqnarray*}
which shows that $F$ is a non constant rational function and the resonance
conditions of definition~\ref{def:res} are fulfilled only on a union of
codimension 1 submanifolds.
Hence it is enough to choose a normal perturbation with $\lambda$ small,
$\lambda '(0)=\lambda '(\pi)=0$ and outside this set to obtain a curve
$\beta$, close to $\alpha$, with a non-resonant elliptic diameter. 
 \ep\medbreak

\begin{prop}\label{prop:nres}
If $0$ is a non-resonant elliptic diameter for a strictly convex $C^k$ curve
$\alpha$, $k\geq 3$, then there exists $\epsilon_2$ such that every normal
perturbation $\beta$, with $|\lambda|_2<\epsilon_2$, has a non-resonant
elliptic diameter.
\end{prop}

\proc{Proof:} : As $0$ is a non-resonant elliptic diameter for $\alpha$ then
\begin{eqnarray*}
\cos\gamma(\alpha) &=& 
2 \, \frac{(l_\alpha(0)-R_\alpha(0)\,(l_\alpha(0)-R_\alpha(\pi))}
{R_\alpha(0)R_\alpha(\pi)}-1 = \\
&=& 2 \, \frac{(L_0-R_0)(L_0-R_\pi)}{R_0R_\pi}-1 \ \neq \  
0, 1, -1, \cos\frac{2}{3}\pi,\cos\frac{4}{3}\pi
\end{eqnarray*}
For $\epsilon_1$ as in Proposition \ref{prop:aberto}, let 
$\beta(\varphi)=\alpha(\varphi)+\lambda(\varphi)\eta(\varphi)$, with
$|\lambda|_2<\epsilon_1$. Then $\beta$ has an elliptic diameter at
$\varphi_1\in (-\delta, \delta)$ with eigenvalues determined by
$$
\cos\gamma(\beta)= 2 \, \frac{(l_\beta(\varphi_1)-R_\beta(\varphi_1))
(l_\beta(\varphi_1)-R_\beta(h^{-1}(h(\varphi_1)+\pi))}
{R_\beta(\varphi_1)R_\beta(h^{-1}(h(\varphi_1)+\pi))}-1.$$
Since $\alpha$ and $\beta$ are $C^k$, $k\geq 3$, then $R_\alpha$ and $R_\beta$
are $C^1$ close. This, together with lemmas \ref{lem:h-1} and \ref{lem:lalfa},
gives the desired result. 
\ep\medbreak

\section{Elliptic Islands}

In order to seek after islands, let us change coordinates from 
$(\varphi,\theta)$ to $(s,p)$, where $s$ is the arclenght parameter for the 
curve $\alpha$ and $p=\sin \theta$.
The advantage of those new coordinates is that the billiard map $T$ preserves
the area $ds\wedge dp$.

Let us now suppose that, for $k\geq 5$, the $C^k$ curve
$\alpha$ has a non resonant elliptic diameter at $s=0$, so that $(0,0)$ is a
non resonant elliptic fixed point of $T^2$.

By means of a linear area preserving coordinate change, and
complexification, $T^2(s,p)$ can be brought into the
form $ e^{i\gamma}z + g(z,\overline z)$ where $z,\overline z$ are complex and
$g$ vanishes with its first derivative at $z=0$. 
This map has then the Birkhoff Normal Form
$e^{i(\gamma+\tau_1(\alpha)|\zeta|^2)}\zeta+h(\zeta,\overline\zeta)$,
where $h(\zeta,\overline\zeta)={\cal O} (|\zeta|^4)$.  Moreover, if the first
Birkhoff coefficient $\tau_1(\alpha)\neq 0$ then, by Moser's Twist Theorem
\cite{kn:mos}, $\zeta=0$ is a stable fixed point, i.e., $(0,0)$ has
$T^2$-invariant curves surrounding it which means that
there is an elliptic island of positive measure. 

Hence, by assuming $\tau_1(\alpha)=0$, the goal is show that small
$C^2$-normal perturbations of $\alpha$ may have non-zero Birkhoff coefficients.

As Moeckel showed in \cite{kn:moe}, $\tau_1$ can be calculated
using complex area preserving coordinates diagonalizing
the linear part of the map, 
such that the third jet of $T^2$ takes the form:
$$
Z=e^{i\gamma}(z+c_{20}z^2+c_{11}z\overline z+c_{02}\overline
z^2+c_{30}z^3+c_{21}z^2\overline z+c_{12}z\overline z^2+c_{03}\overline z^3)
$$
then the Birkhoff coefficient is given by
$$\tau_1=\frac{1}{i}\left (c_{21}+2|c_{20}|^2\left
(\frac{2e^{i\gamma}+1}{e^{i\gamma}-1}+\frac{1}{e^{3i\gamma}-1}\right )
\right) \ .
$$

The first step is, then, to determine the third jet of $T^2$ at $(0,0)$.
Although we are not able to write $T$ explicitly,  we have that (see, for
instance, \cite{kn:str}): \begin{eqnarray*}
\frac{\partial s_1}{\partial s_0}&=&\frac{l(s_0,p_0)-R(s_0)\cos\theta(p_0)}
{R(s_0)\cos\theta(p_1)} 
\\
\frac{\partial s_1}{\partial p_0}&=&\frac{l(s_0,p_0)}
{\cos\theta(p_0)\cos\theta(p_1)} 
\\
\frac{\partial p_1}{\partial
s_0}&=&\frac{l(s_0,p_0)-R(s_0)\cos\theta(p_0)-R(s_1)\cos\theta(p_1)}
{R(s_0)R(s_1)}
\\
\frac{\partial p_1}{\partial p_0}&=&\frac{l(s_0,p_0)-R(s_1)\cos\theta(p_1)}
{R(s_1)\cos\theta(p_0)} 
\end{eqnarray*}
where $(s_1,p_1) = T(s_0,p_0)$. 

Using the Chain Rule and with the aid of the software Maple \footnote{the
program may be found at www.mat.ufmg.br/$\sim$syok}, we have
calculated the coefficients of the third-jet $J_3 \,T^2_{(0,0)}(s_0,p_0)$ and
obtained the above coefficients $c_{ij}$. They depend on $L_0$, $R_0=R(0)$,
$R_{\pi}=R(s_1)$, $\frac{dR}{ds}(0)$, $\frac{dR}{ds}(s_1)$,
$\frac{d^2R}{ds^2}(0)$, $\frac{d^2R}{ds^2}(s_1)$, where $(s_1,0)=T(0,0)$.
Moreover, $c_{21}$ depends
linearly on $\frac{d^2R}{ds^2}(0)$ and $\frac{d^2R}{ds^2}(s_1)$ and $c_{20}$
does not depend on those second derivatives.

Using these calculations we can prove:
\begin{prop}\label{prop:tau}
Let $\alpha$ be a $C^k$ strictly convex curve, $k\geq 5$, with a non-resonant
elliptic diameter at $s=0$, and such that $\tau_1(\alpha)=0$. Then there is a $C^k$ 
strictly convex curve $\beta$,
$C^2$-close to $\alpha$, such that $\beta$ has a non-resonant elliptic
diameter, with $\tau_1(\beta)\neq 0$.
\end{prop}

\proc{Proof:} : Let $\beta(\varphi)=\alpha(\varphi)+\lambda(\varphi)\eta(\varphi)$
be a normal perturbation of $\alpha$, with a third order contact with $\alpha$
at $\varphi(0)=0$ and $\varphi(s_1)=\pi$. This means that
$\lambda(0)=\lambda(\pi)=\lambda
'(0)=\lambda '(\pi)=\lambda ''(0)=\lambda ''(\pi)=\lambda '''(0)=\lambda
'''(\pi)=0$ and then, $\beta$ has a non-resonant elliptic
diameter at $\varphi=0$

Let $\sigma$ be the arclenght parameter for $\beta$ and
$(\sigma_1,0)=T^2_{\beta}(0,0)$. The third order contact implies that the
elliptic diameters coincide, with the same length $L_0$ and
 $$R_{\beta}(0)=R_{\alpha}(0)=R_0,\,\,\,
R_{\beta}(\sigma_1)=R_{\alpha}(s_1)=R_{\pi},$$ 
$$\frac{dR_{\beta}}{d\sigma}(0)=\frac{dR_{\alpha}}{ds}(0),\,\,\,
\frac{dR_{\beta}}{d\sigma}(\sigma_1)=\frac{dR_{\alpha}}{ds}(s_1),$$
$$\frac{dR^2_{\beta}}{d\sigma^2}(0)=\frac{dR^2_{\alpha}}{ds^2}(0)
-\frac{d^4\lambda}{d\varphi^4}(0),\,\,\,
\frac{dR^2_{\beta}}{d\sigma^2}(\sigma_1)=\frac{dR^2_{\alpha}}{ds^2}(s_1)
-\frac{d^4\lambda}{d\varphi^4}(\pi).$$

Then $\gamma(\beta)=\gamma(\alpha)$, $c_{20}(\beta)=c_{20}(\alpha)$ and
$$c_{21}(\beta)=-\frac{i}{16}\left[\frac{L_0(L_0-R_\pi)}{(L_0-R_0)(L_0-R_0-R_\pi)}\right]
\frac{d^4\lambda}{d\varphi^4}(0)+ib(L_0,R_0,R_\pi)\frac{d^4\lambda}{d\varphi^4}(\pi)+ 
c_{21}(\alpha),$$
where $b$ is a $C^{k-4}$ function.

Hence
$$\tau_1(\beta)=-\frac{1}{16}\left[\frac{L_0(L_0-R_\pi)}{(L_0-R_0)(L_0-R_0-R_\pi)}\right]
\frac{d^4\lambda}{d\varphi^4}(0)+b(L_0,R_0,R_\pi)\frac{d^4\lambda}{d\varphi^4}(\pi)+  
\tau_1(\alpha).$$
As $\sigma=0$ is a non-resonant elliptic diameter for $\alpha$ and $\beta$,
then $L_0\neq 0$, $ L_0-R_0\neq 0$, $L_0-R_\pi \neq 0$ and $L_0-R_0-R_\pi\neq
0$. Therefore, by choosing $\frac{d^4\lambda}{d\varphi^4}(0)\neq 0$ and 
$\frac{d^4\lambda}{d\varphi^4}(\pi)=0$ we obtain a perturbation $\lambda$ such
that $\tau_1(\beta)\neq 0$.

\section{Density of Elliptic Islands}

Let $\alpha$ be a closed, regular, strictly convex plane curve and $T_\alpha$ its
associated billiard map. Planar rigid motions do not change the
geometrical features of $\alpha$ and then do not change the dynamical
characteristics of $T_\alpha$. This defines an equivalence relation on the set
of curves. 

Let ${\cal C}$ be the set of equivalent classes, denoted by $[\alpha]$, of
closed, regular, strictly convex plane curves $\alpha$, that are $C^5$ when
parametrized by the angle $\varphi$. 

Given  a representative $\alpha$ of an equivalent class $[\alpha]\in{\cal
C}$, the normal bundle $(\alpha(\varphi),\eta(\varphi))$ is then $C^5$.
Given $\epsilon>0$, let us consider a tubular neighbourhood 
$
N_\epsilon(\alpha)=
\{\alpha(\varphi)+\lambda\eta(\varphi),0\le \varphi < 2\pi,
-\epsilon<\lambda<\epsilon
\} \simeq S^1\times (-\epsilon,\epsilon)$ and the canonical projection  
$\Pi:S^1\times (-\epsilon,\epsilon)\to S^1$ that projects $N_\epsilon(\alpha)$
onto the image of $\alpha$.

\begin{definition}
$[\beta]\in{\cal C}$ is $\epsilon$-close to $[\alpha]\in{\cal C}$ if:
\begin{enumerate}
\item there exists $\beta\in [\beta]$ such that the image of $\beta$ is in
$N_\epsilon(\alpha)$, for $\alpha\in[\alpha]$.
\item the restriction of $\Pi$ to the image of $\beta$ is a diffeomorphism.
\end{enumerate}
\end{definition}
As a consequence, $\beta$ can be written as 
$\beta(\varphi)=\alpha(\varphi)+\lambda(\varphi)\eta(\varphi)$, with
$\lambda$ periodic.

\begin{definition}
$[\beta]\in{\cal C}$ is $\epsilon$-$C^2$-close to $[\alpha]\in{\cal C}$ if
$[\beta]$ is $\epsilon$-close to $[\alpha]$ and $|\lambda|_2<\epsilon$.
\end{definition}

In this context, Propositions \ref{prop:aberto} and \ref{prop:nres} can be
rewritten as:
\begin{prop}
Let ${\cal E}\subset{\cal C}$ be the subset of equivalent classes of curves
with an elliptic diameter. Then ${\cal E}$ is open in ${\cal C}$.
 \end{prop}
\begin{prop}
Let ${\cal R}\subset{\cal E}$ be the subset of equivalent classes of curves
with a non resonant diameter. Then ${\cal R}$ is open in ${\cal C}$.
 \end{prop}

Propositions \ref{prop:reso} and \ref{prop:tau} leads to:

\begin{theorem}
Any curve on an equivalent class of ${\cal E}$ can be approximated by curves
such that the associated billiard map has an elliptic island of positive
measure.  
\end{theorem}

\proc{Proof:} : Let $\alpha$ be a representative of $[\alpha]\in{\cal E}$. So,
$\alpha$ has an elliptic diameter at, we can suppose, $\varphi=0$. There are
three possible cases: 
\begin{enumerate}
\item $0$ is a non resonant diameter with $\tau_1(\alpha)\neq 0$.
\item $0$ is a non resonant diameter with $\tau_1(\alpha)=0$.
\item $0$ is a resonant diameter.
\end{enumerate}
In the first case, there is nothing to prove.

If $0$ is a non resonant diameter with $\tau_1(\alpha)=0$ then, by Proposition
\ref{prop:tau}, $\alpha$ can be approximated by a sequence $\beta_n$ such that
for each $n$, $\beta_n$ has a non resonant diameter with $\tau_1(\beta_n)\neq
0$.

If  $0$ is a resonant diameter, by Proposition \ref{prop:reso}, $\alpha$ can
be approximated by a sequence $\beta_n$ such that for each $n$, $\beta_n$ has a
non resonant diameter. The same reasoning used on the second case leads to
the result. 
\ep\medbreak

\section{Conclusions and Comments}

Although we can prove the density of elliptic islands surrounding orbits
associated to diameters, many questions remain open.

The existence
of islands is not equivalent to a non-vanishing first Birkhoff coefficient.
For instance, any Birkhoff coefficient different from zero assures,
for a sufficiently differentiable curve and avoiding higher order resonances,
the existence of such islands. However one may ask if there is a strictly
convex $C^\infty$ curve with an elliptic diameter such that all the Birkhoff
coefficients vanishes. If so, how this fact translates into the dynamics of the
associated billiard map?

It may also happen that all the isolated diameters of a given curve are
hyperbolic or parabolic. This is the case, for instance, of an ellipse with
half-axes 1 and $\sqrt 2$, which has only two diameters: the largest axis,
which is hyperbolic, and the smallest one, which is parabolic.
However, due to the
integrability of the system, the parabolic orbit has an island.
Kozlov, in \cite{kn:koz} constructs a strictly convex curve, as
differentiable as desired, that has exactly three  diameters, all of them
hyperbolic.

\begin{figure}[ht]
\centering
{\includegraphics[bb = 0 0 340 220]{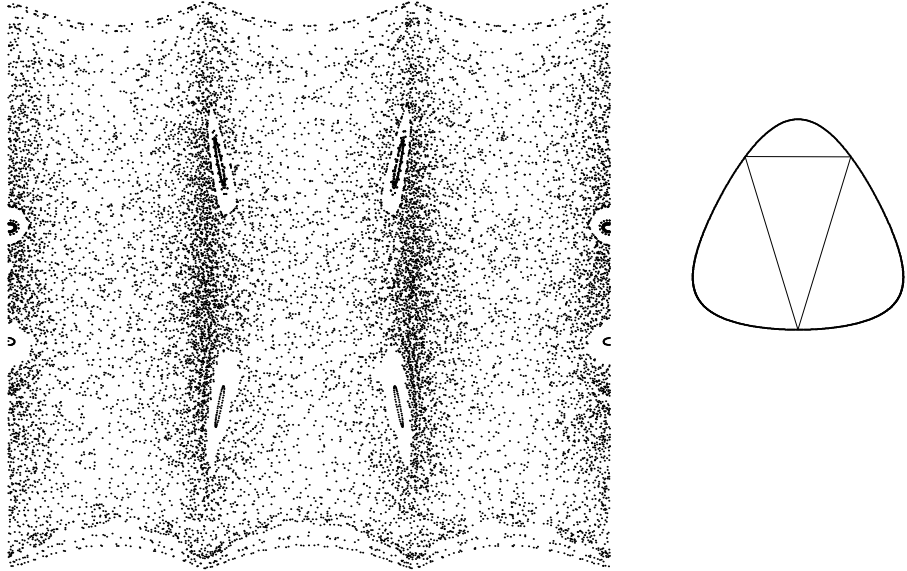}}
\caption{The curve $x(t)=\cos t$, $y(t)=3/(2-\sin t)$
with a period 3 trajectory and the phase space of the associated billiard map.}
\label{fig:curva}
\end{figure}  

Even when all the diameters are hyperbolic, the examples suggest that a higher
period trajectory is elliptic, if the curve is strictly convex. This is the
case for the billiard associated to $x(t)=\cos t$, $y(t)=3/(2-\sin t)$. It has 
exactly two diameters, both hyperbolic. The numerical investigation of the
phase-space of its associated billiard map shows the existence of a pair of
3-periodic orbits, corresponding to a same triangular trajectory followed
in both senses, and its islands (see figure~\ref{fig:curva}).
One may wonder if there is always an elliptic periodic orbit for any
sufficiently  differentiable and strictly convex curve. 

If there are elliptic periodic orbits, can they have any period? 
In the example above only elliptic orbits of period 3 or multiples of 3,
belonging to the same island, are visible.

Anyway, if there are elliptic orbits of higher period and they are
non-resonant, we believe that our method of normal perturbations with a
third order contact at the points of reflection of the periodic trajectory
with the boundary will give the results we obtained for the 2-periodic
orbit.

\end{document}